\magnification=\magstep1
\input amstex
\documentstyle{amsppt}

\font\tencyr=wncyr10 
\font\sevencyr=wncyr7 
\font\fivecyr=wncyr5 
\newfam\cyrfam \textfont\cyrfam=\tencyr \scriptfont\cyrfam=\sevencyr
\scriptscriptfont\cyrfam=\fivecyr
\define\hexfam#1{\ifcase\number#1
  0\or 1\or 2\or 3\or 4\or 5\or 6\or 7 \or
  8\or 9\or A\or B\or C\or D\or E\or F\fi}
\mathchardef\Sha="0\hexfam\cyrfam 58

\define\Primes{\frak{Primes}}

\define\defeq{\overset{\text{def}}\to=}
\define\ab{\operatorname{ab}}
\define\pr{\operatorname{pr}}
\define\Gal{\operatorname{Gal}}
\def \isom {\overset \sim \to \rightarrow}

\define\Spec{\operatorname{Spec}}
\define\id{\operatorname{id}}

\define\Ker{\operatorname{Ker}}
\def \res{\operatorname {res}}

\def\char{\operatorname{char}}
\def\Card{\operatorname{Card}}

\def\sol{\operatorname{sol}}
\def \Res{\operatorname{Res}}
\def \all{\operatorname{all}}

\def \Aut{\operatorname{Aut}}

\def \sep{\operatorname{sep}}

\def \and{\operatorname{and}}

\def \Ind{\operatorname{Ind}}
\def \tame{\operatorname{tame}}

\def \et{\operatorname{et}}

\NoRunningHeads
\NoBlackBoxes
\topmatter

\title
A local-global principle for torsors under geometric prosolvable fundamental groups II
\endtitle

\author
Mohamed Sa\"\i di
\endauthor

\abstract We prove a {\it local-global principle for torsors under the prosolvable geometric fundamental group of
an {\bf affine} curve over a number field}. 
\endabstract
\toc
\subhead
\S 0. Introduction
\endsubhead

\subhead
\S 1. Prosolvable geometric tame fundamental groups
\endsubhead

\subhead
\S 2. Geometrically prosolvable arithmetic fundamental groups
\endsubhead

\subhead
\S 3. Proof of Theorem B
\endsubhead

\endtoc

\endtopmatter

\document

\subhead
\S 0. Introduction
\endsubhead
Let $k$ be a characteristic $0$ field and $X\to \Spec k$ a separated, smooth, and geometrically connected 
{\bf curve} over $k$. Let $\eta$ be a geometric point of $X$ with values in its generic point.
Thus, $\eta$ determines an algebraic closure $\bar k$ of $k$, and a geometric point $\overline \eta$
of $X_{\bar k}\defeq X\times _{\Spec k}\Spec \bar k$.
There exists a canonical exact sequence of profinite groups (cf. [Grothendieck], Expos\'e IX, Th\'eor\`eme 6.1)
$$1\to \pi_1(X_{\bar k},\overline \eta)\to \pi_1(X, \eta) @>>> G_k\to 1.$$
Here, $\pi_1(X, \eta)$ denotes the arithmetic \'etale fundamental group of $X$ with base
point $\eta$, $\pi_1(X_{\bar k},\overline \eta)$ denotes the \'etale fundamental group of $X_{\bar k}$ with base
point $\overline \eta$, and $G_k\defeq \Gal (\overline k/k)$ denotes the absolute Galois group of $k$.
Write
$$\Delta\defeq \pi_1(X_{\bar k},\overline \eta)\ \ \ \ \ \ \ \ \ \  \text {and}\ \ \ \ \ \ \ \ \ \ \Pi\defeq \pi_1(X, \eta).$$
We have the above natural exact sequence 
$$ 1\to \Delta \to \Pi\to G_k\to 1.\tag 0.1$$

Suppose that the sequence $(0.1)$ splits
[for example assume that $X(k)\neq \emptyset$]. Let 
$$s:G_k\to \Pi$$ 
be a section of the projection $\Pi \twoheadrightarrow G_k$. We view $\Delta$ as a $G_k$-group
via the conjugation action of $s(G_k)$.

Assume that $k$ is a {\bf number field}; i.e., $k$ is a finite extension of $\Bbb Q$. 
Let $v$ be a prime of $k$, $k_v$ the completion of $k$ at $v$, and $G_{k_v}\subset G_k$ a decomposition group associated to $v$. 
Thus, $G_{k_v}$ is only defined up to conjugation. We view $\Delta$ as a $G_{k_v}$-group
via the conjugation action of $s(G_{k_v})$. For each prime $v$ of $k$ we have a natural restriction map (of pointed non-abelian cohomology sets) 
$$\Res_v:H^1(G_k,\Delta)\to H^1(G_{k_v},\Delta),$$
and a natural map
$$\prod _{\all v} \Res _v: H^1(G_k,\Delta)\to \prod _{\all v} H^1(G_{k_v},\Delta),$$
where the product is over all primes $v$ of $k$. The main  problem we are concerned with in this paper is the following.

\definition {Question A} Is the map 
$$\prod _{\all v} \Res _v: H^1(G_k,\Delta)\to \prod _{\all v} H^1(G_{k_v},\Delta)$$
{\bf injective}?
\enddefinition

As explained in [Sa\"\i di] $\S0$, the above question is related to the Grothendieck anabelian section conjecture.
Let $\Delta^{\sol}$ be the maximal {\bf prosolvable} quotient of $\Delta$, which is a characteristic quotient.
The above $G_k$ (resp. $G_{k_v}$)-group structure on $\Delta$ induces naturally a  $G_k$ (resp. $G_{k_v}$)-group structure on
$\Delta^{\sol}$. Let $\Primes_k$ be the set of primes of $k$ and $S\subseteq \Primes _k$ a non-empty subset, we have as above a natural restriction map

$$\prod _{v\in S} \Res _v^{\sol}: H^1(G_k,\Delta^{\sol})\to \prod _{v\in S} H^1(G_{k_v},\Delta^{\sol}).$$

In [Sa\"\i di] we proved the following result.

\proclaim {Theorem A} Assume $k$ is a {\bf number field}, $X$ is {\bf proper},
and $S\subseteq \Primes _k$ is a set of primes of $k$ of {\bf density $1$}.
Then the map
$$\prod _{v\in S} \Res _v^{\sol}: H^1(G_k,\Delta^{\sol})\to \prod _{v\in S} H^1(G_{k_v},\Delta^{\sol})$$
is {\bf injective}.
\endproclaim

In this note we generalise Theorem A by removing the assumption therein that $X$ is proper. Our main result is the following.
\proclaim {Theorem B} Assume $k$ is a {\bf number field}, $X$ is {\bf affine}, and $S\subseteq \Primes _k$ is a set of primes of $k$ of {\bf density $1$}.
Then the map
$$\prod _{v\in S} \Res _v^{\sol}: H^1(G_k,\Delta^{\sol})\to \prod _{v\in S} H^1(G_{k_v},\Delta^{\sol})$$
is {\bf injective}.
\endproclaim

Our proof of Theorem B relies on a devissage argument, and a careful analysis of the structure of the geometric prosolvable 
(resp. geometrically prosolvable arithmetic) (tame) fundamental group of an affine curve which is established in $\S1$ (resp. $\S2$). The results established in $\S1$ and $\S2$ 
may be of interest independently of the question discussed in this paper. 
The abelian analog of Theorem B is a consequence of results of Serre (cf. proof of Proposition 3.2). Theorem B is proved in $\S3$. 

\smallskip
{\bf Acknowledgment.} I am very grateful to the referee for his/her careful reading of the paper and valuable comments.

\bigskip
\subhead
Notations
\endsubhead
The following notations will be used throughout this paper (unless we specify otherwise).
For a profinite group $H$:

\smallskip
$\bullet$ we denote by $\overline {[H,H]}$ the closed subgroup of $H$ which is 
topologically generated by the commutator subgroup of $H$.

$\bullet$ we denote by $H^{\ab}\defeq  {H}/{\overline {[H,H]}}$ the maximal abelian quotient of $H$.

$\bullet$ we denote by $H^{\sol}$ the maximal {\bf prosolvable} quotient of $H$.

\smallskip
$\bullet$ For an exact sequence $1\to H'\to H@>\pr>> G\to 1$ of profinite groups we will refer to a continuous homomorphism
$s:G\to H$ such that $\pr \circ s=\id_G$ as a group-theoretic {\bf section}, or simply section, of the natural projection $\pr:H\twoheadrightarrow G$.

\subhead
\S1. Prosolvable geometric tame fundamental groups
\endsubhead
Let $\ell$ be an {\bf algebraically closed} field of characteristic $p\ge 0$, $V\to \Spec \ell$ a separated, smooth, and connected 
$\ell$-{\bf curve}, $Z$ the smooth compactification of $V$, and $Z\setminus V$ the complement of $V$ in $Z$ which is either empty
if $Z=V$ or otherwise $V$ is affine and $Z\setminus V$ consists of finitely many closed points of $Z$.
Let $r\defeq \Card (Z\setminus V)$ be the cardinality of $Z\setminus V$, thus $r=0$ if $Z=V$.  
Let $\xi$ be a geometric point of $V$ with values in its generic point, and 
$$\Delta\defeq \pi_1^{\tame}(V,\xi)$$ 
the {\bf tame fundamental group} of $V$ with base point $\xi$.
The geometric point $\xi$ determines a separable closure $K^{\sep}$ of the function field $K$ of $Z$.
Thus, $\Delta\defeq \Gal (K^{\tame}/K)$ where $K^{\tame}$ is the maximal sub-extension of $K^{\sep}/K$ which is 
\'etale above $V$ and tamely ramified over the (discrete) valuations of $K$ corresponding to closed points of $Z\setminus V$.

Consider the {\bf derived series} of $\Delta$
$$.....\subseteq \Delta (i+1)\subseteq \Delta (i)\subseteq......\subseteq \Delta (1)\subseteq \Delta(0)=\Delta \tag 1.1$$
where 
$$\Delta(i+1)=\overline {[\Delta (i),\Delta(i)]},$$
for $i\ge 0$, is the $i+1$-th derived subgroup which is a characteristic subgroup of $\Delta$. Write
$$\Delta _i\defeq \Delta/\Delta(i).$$
Thus, $\Delta _i$ is the {\bf $i$-th step solvable} quotient of $\Delta$, and $\Delta _1\defeq \Delta^{\ab}$ is the maximal {\bf abelian} quotient of $\Delta$. 
Note that there exists a natural exact sequence
$$ 1\to   \Delta^{i+1}\to \Delta _{i+1}\to \Delta _i\to 1 \tag 1.2$$
where $\Delta^{i+1}$ is the subgroup $\Delta (i)/\Delta (i+1)$ of $\Delta _{i+1}$. In particular, $\Delta ^{i+1}$ is {\bf abelian}.

\proclaim {Lemma 1.1} We have a natural identification $\Delta ^{\sol}  \isom \underset{i\ge 1} \to{\varprojlim}\ \Delta_i$. 
In particular, $\Ker (\Delta \twoheadrightarrow \Delta ^{\sol})=\bigcap _{i\ge 1}\Delta (i)$,
and $\Delta ^{\sol}$ is a characteristic quotient of $\Delta$.
\endproclaim

\demo{Proof} Follows from the various definitions.
\qed
\enddemo

Let $i\ge 1$ be an integer. The profinite group $\Delta _{i}$ is finitely generated as follows from the well-known
finite generation of $\Delta$ which projects onto $\Delta _{i}$ (cf. [Grothendieck], Expos\'e XIII, Corollaire 2.12).  
Let $\{\widehat \Delta _{i}[n]\}_{n\ge 1}$ be a countable system of {\bf characteristic open} subgroups of $\Delta_i$ such that 
$$\widehat \Delta_{i}[n+1]\subseteq \widehat \Delta _{i}[n],\ \ \ \ \widehat \Delta _{i}[1]\defeq \Delta_i,\ \ \ \ \text {and} \ \ \bigcap _{n\ge 1}\widehat \Delta_{i}[n]=\{1\}.$$

Write ${\Delta_i}[n]$ for the inverse image of $\widehat \Delta _i[n]$ in $\Delta$ via the map $\Delta \twoheadrightarrow \Delta_i$. Thus, $\Delta_i[n]\subseteq \Delta$ is an open subgroup corresponding to an \'etale Galois cover $V_{i,n}\to V$ which extends to a tamely ramified Galois cover $Z_{i,n}\to Z$ between proper, smooth, and connected 
$\ell$-curves (the cover $Z_{i,n}\to Z$  is \'etale if $r=0$, i.e., if $V=Z$, or if $i=r=1$). Note that since the \'etale cover $V_{i,n}\to V$ 
is defined via an open subgroup of $\pi_1^{\tame}(V,\xi)$ it is a pointed \'etale cover, and $V_{i,n}$ 
(hence also $Z_{i,n}$) is naturally endowed with a geometric point $\xi_{i,n}$ above $\xi$.

If $i=r=1$, let
$$\widetilde {\Delta}_i[n]\defeq  {\Delta}_i[n]=\pi_1^{\tame}(V_{i,n},\xi_{i,n})$$ 
be the tame fundamental group of $V_{i,n}$ with base point $\xi_{i,n}$. 

If $(i,r)\neq (1,1)$, let 
$$\widetilde {\Delta}_i[n]\defeq \pi_1(Z_{i,n},\xi_{i,n})$$ 
be the \'etale fundamental group of 
$Z_{i,n}$ with geometric point $\xi_{i,n}$, which is a quotient of ${\Delta_i}[n]=\pi_1^{\tame}(V_{i,n},\xi_{i,n})$. 

Let $i\ge 1$. We have, $\forall n\ge 1$, finite morphisms $V_{i,n+1}\to V_{i,n}$, and
$Z_{i,n+1}\to Z_{i,n}$, which induce continuous homomorphisms
$\widetilde {\Delta} _i[n+1]\to \widetilde {\Delta} _i[n]$
and
$\widetilde {\Delta} _i[n+1]^{\ab}\to \widetilde {\Delta} _i[n]^{\ab}$. Thus, we have a projective system $\{\widetilde {\Delta} _i[n]^{\ab} \}_{n\ge 1}$.

Our main result in this section is the following description of the structure of $\Delta ^{i+1}$; this description in the case  
$(i,r)\neq (1,1)$ is in terms of the Tate modules
of the jacobians of the $\{Z_{i,n}\}_{n\ge 1}$. Proposition 1.2 
plays a key role in the proof of Theorem B and may be of interest independently of the question discussed in this paper.

\proclaim {Proposition 1.2} Assume $i\ge 1$. We have a natural identification 
$$\Delta ^{i+1}  \isom \underset{n\ge 1} \to{\varprojlim}\ \widetilde \Delta_i[n] ^{\ab}.$$
\endproclaim

\demo {Proof} If $r=0$, i.e., $V=Z$ is proper, or if $i=r=1$, the assertion follows easily from the various Definitions. Observe that in the case $i=r=1$ the natural projection 
$\Delta_1=\pi_1^{\tame}(V,\xi)^{\ab}\twoheadrightarrow \pi_1(Z, \xi)^{\ab}$ is an isomorphism.

Next, we assume $V=Z\setminus \{x_1,\ldots,x_r\}$ is affine and $V$ is the complement in $Z$
of a finite set $\{x_1,\ldots,x_r\}$ of closed points with $\max (i,r)>1$.
Let $G$ be a finite quotient of $\Delta_{i+1}$ which inserts in 
the following commutative diagram
$$
\CD
1 @>>> \Delta ^{i+1}  @>>>  \Delta _{i+1} @>>> \Delta _i @>>> 1\\
@. @VVV     @VVV   @VVV\\
1 @>>> A  @>>> G @>>> G' @>>> 1\\
\endCD
$$
where the vertical maps are surjective.
The quotient $G$ corresponds to a finite \'etale Galois cover $V_1\to V$ which extends to a tamely ramified
finite Galois cover $Z_1\to Z$ with Galois group $G$. The cover $Z_1\to Z$ factorizes as
$Z_1\to Z_1'\to Z$, where $Z_1'\to Z$ is the sub-cover with Galois group $G'$, and $Z_1\to Z_1'$ is a Galois abelian cover with group $A$. 
For $s\in \{1,\ldots,r\}$, let $I_{x_s}\subset G$ be an inertia subgroup associated to $x_s$. Thus, $I_{x_s}$ is only defined up to conjugation.
Moreover, $I_{x_s}$ is cyclic of order $e_s\ge 1$ with $\text {gcd}(e_s,p)=1$ as the ramification is tame. 
The following claim follows immediately from the well-known structure of $\Delta$ (cf. [Grothendieck], Expos\'e XIII, Corollaire 2.12). 

\proclaim {Claim 1.2.2} There exists a finite quotient $\Delta _{i}\twoheadrightarrow P$ of $\Delta_{i}$ 
corresponding  to a finite tamely ramified Galois cover $Z_2\to Z$ which is \'etale above $V$
such that the followings hold. For $s\in \{1,\ldots,r\}$, if $I'_{x_s}\subseteq P$ is an inertia subgroup associated to $x_s$, 
then $I'_{x_s}$ is {\it cyclic} of order $f_s=e_sh_s$ a multiple of $e_s$ with $\text {gcd}(f_s,p)=1$.
\endproclaim

Note that Claim 1.2.2, which holds under the assumption $\max (i,r)>1$, is false if $i=r=1$. Indeed, a finite abelian tamely ramified cover $C'\to C$ 
between proper and connected smooth $\ell$-curves which is \'etale above $C\setminus \{x\}$, where $x\in C$ is a closed point, is necessarily \'etale above $C$.

Next, let $K_1\defeq K_{Z_1}$ (resp. $K_2\defeq K_{Z_2}$) be the function field of $Z_1$ (resp. $Z_2$), 
$L\defeq K_1.K_2$ the compositum of $K_1$ and $K_2$ in $K^{\sep}$, and $\widetilde Z$ the normalisation of $Z$ in $L$. Thus, $\widetilde Z\to Z$ is a 
tamely ramified Galois cover with Galois group $H\subseteq G\times P$ which is \'etale above $V$.

\proclaim {Lemma 1.2.3}  The quotient $\Delta \twoheadrightarrow H$ factorizes as   $\Delta \twoheadrightarrow \Delta _{i+1}\twoheadrightarrow H$. 
\endproclaim

\demo{Proof of Lemma 1.2.3} This follows from the facts that $\widetilde Z\to Z$ is \'etale above $V$,
tamely ramified, $H\subseteq G\times P$, and $G\times P$ is $(i+1)$-th step solvable.
\qed
\enddemo

Let $H'$ be the image of $H$ in $G'\times P$. We have a commutative diagram of exact sequences where the vertical maps are natural inclusions.
$$
\CD
1  @>>> A_H\defeq H\cap \lgroup A\times \{1\} \rgroup @>>> H @>>> H' @>>>1\\
@. @VVV  @VVV  @VVV  \\
1 @>>> A\times \{1\}  @>>> G\times P @>>> G'\times P @>>> 1\\
\endCD
$$

\proclaim {Lemma 1.2.4}
The group $H'$ is a quotient of $\Delta _{i}$. Moreover, the cover $\widetilde Z\to Z$ factorizes as $\widetilde Z\to \widetilde Z'\to Z$
where $\widetilde Z'\to Z$ is \'etale above $V$, tamely ramified, and Galois with Galois group $H'$, and $\widetilde Z\to \widetilde Z'$ is an {\bf abelian \'etale} cover with Galois group $A_H$.
\endproclaim

\demo{Proof of Lemma 1.2.4} The first assertion, as well as the assertion concerning the cover $\widetilde Z'\to Z$, follow from the various definitions. 
We prove the last assertion regarding the cover $\widetilde Z\to \widetilde Z'$. Recall the factorisation
$Z_1\to Z_1'\to Z$ where $Z_1\to Z_1'$ is Galois with abelian Galois group $A$ and 
$Z_1'\to Z$ is Galois with group $G'$. Let $\widetilde Z'$ be the normalisation of $Z$ in the compositum of the function fields of  $Z_1'$ and $Z_2$ (in $K^{\sep}$).
Thus, $\widetilde Z'\to Z$ is a Galois cover with Galois group $H'$ and
we have the following commutative diagram
$$
\CD
\widetilde Z @>>> Z_1\\
@VVV     @VVV\\
\widetilde Z' @>>>  Z_1'\\
@VVV   @VVV\\
Z_2   @>>> Z\\
\endCD
$$
of finite Galois covers.
The ramification index in the Galois cover $Z_2 \to Z$ 
above a branched (closed) point $x_s\in Z$ is divisible by the ramification index above $x_s$ in the Galois cover $Z_1\to Z$ (cf. the condition in Claim 1.2.2 
that $f_s$ is divisible by $e_s$). The fact that the cover $\widetilde Z\to Z_2$, and a fortiori $\widetilde Z\to \widetilde Z'$ which is abelian with Galois group $A_H$, is \'etale follows from Abhyankar's Lemma (cf. [Grothendieck], Expos\'e X, Lemma 3.6).
\qed
\enddemo

Going back to the proof of Proposition 1.2, the above discussion shows that the finite quotients $\Delta _{i+1}\twoheadrightarrow H$ as in Lemma 1.2.3 form a cofinal system 
of finite quotients of $\Delta _{i+1}$. Thus,
$\Delta _{i+1}\isom \underset{H} \to{\varprojlim}\ H$. Proposition 1.2 follows from the facts that the various $H$ above fit in an exact sequence
$1\to A_H\to H\to H'\to 1$, $\Delta _{i}\isom \underset{H'} \to{\varprojlim}\ H'$, and the above Galois covers $\widetilde Z\to \widetilde Z'$ 
with group $A_H$ are \'etale abelian (cf. Lemma 1.2.4).

This finishes the proof of Proposition 1.2.
\qed
\enddemo

\subhead
\S2. Geometrically prosolvable arithmetic fundamental groups
\endsubhead
In this section we use the notations in $\S0$: $k$ is a {\bf field of characteristic $0$}, and $X\to \Spec k$ is a separated, smooth, and geometrically connected (not necessarily proper)
{\bf curve}. 
Consider the {\bf derived series} of $\Delta$ [recall $\Delta=\pi_1(X_{\bar k},\bar \eta)$]
$$.....\subseteq \Delta (i+1)\subseteq \Delta (i)\subseteq......\subseteq \Delta (1)\subseteq \Delta(0)=\Delta \tag 2.1$$
where
$$\Delta(i+1)=\overline {[\Delta (i),\Delta(i)]}, \ \ \ \ \forall i\ge 0.$$

For $i\ge 0$ write
$$\Delta _i\defeq \Delta/\Delta(i).$$
Thus, $\Delta _i$ is the {\bf $i$-th step solvable} quotient of $\Delta$, and $\Delta _1\defeq \Delta^{\ab}$ is the maximal abelian quotient of $\Delta$. 
There exists a natural exact sequence
$$ 1\to   \Delta^{i+1}\to \Delta _{i+1}\to \Delta _i\to 1 \tag 2.2$$
where $\Delta^{i+1}$ is the subgroup $\Delta (i)/\Delta (i+1)$ of $\Delta _{i+1}$ (cf. $\S1$).

Recall $\Pi\defeq \pi_1(X,\eta)$ [cf. exact sequence (0.1)]. Denote 
$$\Pi_i\defeq \Pi/\Delta (i)$$
which inserts in the following exact sequence 
$$1\to \Delta_i\to \Pi_i\to G_k\to 1.$$
We will refer to $\Pi_i$ as the {\bf geometrically $i$-th step solvable} fundamental group of $X$.
We have natural commutative diagrams
of exact sequences

$$
\CD
1  @>>>  \Delta   @>>> \Pi  @>>> G_k  @>>> 1  \\
@.   @VVV       @VVV       @| \\
1  @>>> \Delta _i  @>>>  \Pi_i   @>>>  G_k  @>>> 1 \\
\endCD
\tag 2.3
$$

and

$$
\CD
@.  1  @.  1  \\
@. @VVV    @VVV\\
@. \Delta ^{i+1}   @=  \Delta ^{i+1}  \\
@.  @VVV  @VVV\\
1   @>>>   \Delta _{i+1} @>>>  \Pi_{i+1}  @>>> G_k @>>> 1\\
@.       @VVV            @VVV             @| \\
1 @>>> \Delta _i   @>>> \Pi_i  @>>>  G_k @>>> 1\\
@.     @VVV       @VVV        @VVV \\
   @. 1 @. 1@. 1@. \\
\endCD
\tag 2.4
$$   
where the left and middle vertical maps in diagram (2.3) are the natural surjections.

Consider the pushout diagram
$$
\CD
1 @>>> \Delta  @>>>  \Pi @>>> G_k  @>>> 1\\
@. @VVV      @VVV    @|  \\
1 @>>>  \Delta ^{\sol}   @>>> \Pi^{(\sol)}   @>>> G_k @>>> 1\\
\endCD
\tag 2.5
$$
which defines the {\bf geometrically prosolvable} fundamental group $\Pi^{(\sol)}$ of $X$.

\proclaim {Lemma 2.1} We have a natural identification $\Pi ^{(\sol)}\isom \underset{i\ge 1} \to{\varprojlim}\  \Pi_i$. 
\endproclaim

\demo{Proof} Follows from Lemma 1.1 and the various definitions.
\qed
\enddemo

For the rest of this section we will assume that $X$ is {\bf affine} and denote by $Y$ the {\bf smooth compactification} of $X$.
Thus, $Y\setminus X=\{P_1,\ldots,P_m\}$ consists of $m$-distinct closed points of $Y$, $m\ge 1$. The geometric point $\eta$ of $X$ induces geometric points
$\eta$ and $\bar \eta$ of $Y$, and $Y_{\bar k}\defeq Y\times _{\Spec k}\Spec \bar k$; respectively. We have a natural exact sequence of fundamental groups
$$1\to  \Delta^{\et}\defeq \pi_1(Y_{\bar k},\bar \eta) \to \Pi^{\et}\defeq \pi_1(Y,\eta) \to G_k\to 1.$$ 

By pushing this sequence by the natural projection  
$\Delta ^{\et} \twoheadrightarrow \Delta^{\et}_1\defeq \pi_1(Y_{\bar k},\bar \eta)^{\ab}$ we obtain an exact sequence
$$1\to  \Delta_1^{\et} \to \Pi_1^{\et}\to G_k\to 1$$ 
where $\Pi_1^{\et}\defeq \pi_1(Y,\eta)^{(\ab)}$ is the {\bf geometrically abelian} quotient of 
$\pi_1(Y,\eta)$.  We have an exact sequence 
$$1\to I_X\to \Pi_1\to \Pi_1^{\et}\to 1\tag 2.6$$
where $I_X\defeq \Ker (\Pi_1\twoheadrightarrow  \Pi_1^{\et})=\Ker (\Delta_1 \twoheadrightarrow 
 \Delta_1^{\et})$ is the {\bf inertia subgroup} of $\Pi_1$.
Further, we have an exact sequence of $G_k$-modules
$$0\to \hat \Bbb Z(1)@>>> \prod_{i=1}^m \Ind_{k(P_i)}^k \hat \Bbb Z(1)\to I_X\to 0\tag 2.7$$ 
as follows from the well-known structure of $\Delta_1=\pi_1(X_{\overline k},\bar \eta)^{\ab}$ (cf. the discussion in [Sa\"\i di1], $\S0$). 
Here $\hat \Bbb Z\defeq \prod _{l}\Bbb Z_l$ where the product is over all prime numbers
and $(1)$ is a Tate twist.

Next, let $i\ge 1$ be an integer. The profinite group $\Delta _{i}$ is finitely generated as follows from the well-known
finite generation of $\Delta$ which projects onto $\Delta _{i}$ (cf. the discussion after Lemma 1.1 and the references therein).  
Let $\{\widehat \Delta _{i}[n]\}_{n\ge 1}$ be a countable system of {\bf characteristic open} subgroups of $\Delta_i$ such that 
$$\widehat \Delta_{i}[n+1]\subseteq \widehat \Delta _{i}[n],\ \ \ \ \widehat \Delta _{i}[1]\defeq \Delta_i,\ \ \ \ \text {and} \ \ \bigcap _{n\ge 1}\widehat \Delta_{i}[n]=\{1\}.$$

Write $\Delta _{i,n}\defeq \Delta_i/\widehat \Delta_i[n]$. Thus, $\Delta _{i,n}$ is a finite characteristic quotient of $\Delta _i$
which is an i-th step solvable group. We have a pushout diagram of exact sequences

$$
\CD
1@>>> \Delta_i @>>>   \Pi_i  @>>>  G_k @>>> 1\\
@.  @VVV    @VVV    @| \\
1@>>>  \Delta _{i,n}  @>>>  \Pi _{i,n} @>>>  G_k @>>> 1\\
\endCD
\tag 2.8
$$
where the lower sequence defines a (geometrically finite) quotient $\Pi_{i,n}$ of $\Pi_i$.

In the following discussion we {\bf fix} an integer $i\ge 1$. Suppose that the exact sequence $1\to \Delta_i\to \Pi_i\to G_k\to 1$ {\bf splits}. Let
$$s_i:G_k\to \Pi_i$$
be a section of the natural projection $\Pi_i\twoheadrightarrow G_k$
which induces a section
$$s_{i,n}:G_k\to \Pi_{i,n}$$ 
of the projection $\Pi_{i,n}\twoheadrightarrow G_k$ [cf. diagram (2.8)], for each $n\ge 1$.
Write 
$${\widehat \Pi _i}[n] \defeq {\widehat \Pi _i}[n][s_i] \defeq \widehat \Delta _i[n]. s_i (G_k).$$

Thus, ${\widehat \Pi_i}[n]\subseteq \Pi_i$ is an open subgroup which contains the image 
$s_i(G_k)$ of $s_i$. Write 
$${\Pi_i}[n]\defeq {\Pi_i}[n][s_i]$$ 
for the inverse image of $\widehat \Pi _i[n]$
in $\Pi$ [cf. diagram (2.3)]. Thus, $\Pi_i[n]\subseteq \Pi$ is an open subgroup corresponding to an \'etale cover $X_{i,n}\to X$, which extends to a (possibly ramified) cover 
$$Y_{i,n}\to Y$$
between proper and smooth curves with $Y_{i,n}$ geometrically irreducible (since $\Pi_i[n]$ maps onto $G_k$ via the natural projection $\Pi\twoheadrightarrow G_k$ by the very 
definition of $\Pi_i[n]$). 
Note that the finite cover $(Y_{i,n})_{\bar k}\defeq Y_{i,n}\times _{\Spec k}\Spec \overline k\to Y_{\bar k}$ is Galois with
Galois group $\Delta _{i,n}$, and we have a commutative diagram of covers

$$
\CD
(Y_{i,n})_{\bar k}  @>>>   Y_{\bar k} \\
@VVV    @VVV \\
Y_{i,n}  @>>>  Y\\
\endCD
$$
where $(Y_{i,n})_{\bar k}\to Y$ is Galois with Galois group $\Pi_{i,n}$ and $(Y_{i,n})_{\bar k}\to Y_{i,n}$ is Galois with Galois group $s_{i,n}(G_k)$.

Let $r$ be the number of geometric points of $Y\setminus X$, i.e., $r$ is the cardinality of the finite set of closed points of $(Y\setminus X)_{\bar k}\defeq (Y\setminus X) 
\times _{\Spec k}\Spec \bar k$, where we view $Y\setminus X$ as a (reduced) closed sub-scheme of $Y$ (thus $r\ge m$). 

If $(i,r)=(1,1)$, let
$$\widetilde \Delta _i[n]\defeq \pi_1((X_{i,n})_{\bar k},\overline \eta_{i,n})$$ 
where $(X_{i,n})_{\bar k}\defeq X_{i,n}\times _{\Spec k}\Spec \overline k$, and
$$\widetilde \Pi_i[n]\defeq \pi_1(X_{i,n},\eta_{i,n}).$$

If $(i,r)\neq (1,1)$, let
$$\widetilde \Delta _i[n]\defeq \pi_1((Y_{i,n})_{\bar k},\overline \eta_{i,n})$$ 
and
$$\widetilde \Pi_i[n]\defeq \pi_1(Y_{i,n},\eta_{i,n}).$$
We have, $\forall n\ge 1$, finite morphisms $X_{i,n+1}\to X_{i,n}$, and $Y_{i,n+1}\to Y_{i,n}$, 
 which induce a commutative diagram of exact sequences
$$
\CD
1 @>>>  \widetilde \Delta _i[n+1]@>>> \widetilde \Pi_i[n+1] @>>>  G_k @>>> 1\\
@.   @VVV    @VVV   @| \\
1 @>>>  \widetilde \Delta _i[n]@>>> \widetilde \Pi_i[n] @>>>  G_k @>>> 1\\\
\endCD
\tag 2.9
$$

Note that since the \'etale cover $X_{i,n}\to X$ [resp. $(X_{i,n})_{\bar k}\to X_{\bar k}$] 
is defined via an open subgroup of $\pi_1(X,\eta)$ [resp. $\pi_1({X}_{\bar k},\overline \eta)$] it is a pointed \'etale cover and   $X_{i,n}$ [resp. $(X_{i,n})_{\bar k}$], 
hence also $Y_{i,n}$ [resp. $(Y_{i,n})_{\bar k}$)], 
is naturally endowed with a geometric point $\eta_{i,n}$ (resp. $\overline \eta_{i,n}$) above $\eta$ (resp. $\overline \eta$).

For each integer $n\ge 1$, consider the pushout diagram
$$
\CD
1 @>>>  \widetilde \Delta _i[n] @>>> \widetilde \Pi_i[n] @>>>  G_k @>>> 1\\
@.  @VVV    @VVV   @| \\
1 @>>>  \widetilde \Delta _i[n]^{\ab} @>>> \widetilde \Pi_i[n]^{(\ab)} @>>>  G_k @>>> 1\\
\endCD
\tag 2.10
$$
where $\widetilde \Pi_i[n]^{(\ab)}=\pi_1(Y_{i,n},\eta_{i,n})^{(\ab)}$ (resp. $\widetilde \Pi_i[n]^{(\ab)}=\pi_1(X_{i,n},\eta_{i,n})^{(\ab)}$, if $i=r=1$)
is the geometrically abelian fundamental group of $Y_{i,n}$ (resp. $X_{i,n}$). Further, consider the following commutative diagram

$$
\CD
1 @>>>  \Delta^{i+1}   @>>>  \Cal H_i\defeq \Cal H_{i}[s_i] @>>>  G_k  @>>> 1\\
@.    @|     @VVV     @V{s_i}VV \\
1  @>>>   \Delta ^{i+1}  @>>>   \Pi _{i+1}   @>>>  \Pi_i @>>> 1\\
\endCD
\tag 2.11
$$
where the lower exact sequence is the sequence in diagram (2.3) and the right square is {\bf cartesian}. Thus, (the group extension) $\Cal H_i$
is the pullback of (the group extension) $\Pi_{i+1}$ via the section $s_i$.

\proclaim {Proposition 2.2} Assume $i\ge 1$. We have natural identifications 
$$\Delta ^{i+1}  \isom \underset{n\ge 1} \to{\varprojlim}\ \widetilde \Delta_i[n] ^{\ab}\ \ \ \ \ \text{and}\ \ \ \ \ \ \Cal H_i\isom \underset{n\ge 1} \to{\varprojlim}\  \widetilde \Pi_i[n]^{(\ab)}.$$
\endproclaim

\demo{Proof} The first assertion is Proposition 1.2 [note that $\Delta=\pi_1(V,\eta)=\pi_1^{\tame}(V,\eta)$ since $\char(k)=0$]. 
The second assertion follows from the first and the various Definitions.
More precisely, assuming $(i,r)\neq (1,1)$, let $\widetilde Y_{i,n}\to Y_{i,n}$ be the pro-\'etale cover with Galois group  $\widetilde \Pi_i[n]^{(\ab)}=\pi_1(Y_{i,n},\eta_{i,n})^{(\ab)}$.
Then $\widetilde Y_{i,n}\to Y$ is a Galois cover with Galois group $\widetilde \Pi_{i,n}$. 
We have a commutative diagram of morphisms
$$
\CD
\widetilde Y_{i,n} @>>> (Y_{i,n})_{\bar k}  @>>>   Y_{\bar k} \\
@. @VVV    @VVV \\
@. Y_{i,n}  @>>>  Y\\
\endCD
$$
which induces the following commutative diagram of exact sequences (recall the morphism $(Y_{i,n})_{\bar k}\to Y$
is Galois with Galois group $\Pi_{i,n}$)

$$
\CD
@.  1  @.  1  \\
@. @VVV    @VVV\\
@. \widetilde \Delta _i[n]^{\ab}   @=  \widetilde \Delta _{i}[n]^{\ab}  \\
@.  @VVV  @VVV\\
1   @>>>   \widetilde \Delta _{i,n} @>>>  \widetilde \Pi_{i,n}  @>>> G_k @>>> 1\\
@.       @VVV            @VVV             @| \\
1 @>>> \Delta _{i,n}   @>>> \Pi_{i,n}  @>>>  G_k @>>> 1\\
@.     @VVV       @VVV        @VVV \\
   @. 1 @. 1@. 1@. \\
\endCD
$$   
where $\widetilde \Delta _{i,n} \defeq  \Ker (\widetilde \Pi_{i,n}  \twoheadrightarrow G_k)$. 
Recall the section $s_{i,n}:G_k\to \Pi_{i,n}$ of the projection  $\Pi_{i,n}  \twoheadrightarrow G_k$
induced by the section $s_i$ [cf. discussion after diagram (2.8)] and consider the commutative diagram

$$
\CD
1 @>>>  \widetilde \Delta _i[n]^{\ab}    @>>>  H_{i,n}\defeq H_{i,n}[s_{i,n}] @>>>  G_k  @>>> 1\\
@.    @|     @VVV     @V{s_{i,n}}VV \\
1  @>>>  \widetilde \Delta _i[n]^{\ab}    @>>> \widetilde \Pi_{i,n}   @>>>  \Pi_{i,n} @>>> 1\\
\endCD
$$
where the right square is cartesian. Thus, (the group extension) $H_{i,n}$
is the pullback of (the group extension) $\widetilde \Pi_{i,n}$ via the section $s_{i,n}$. It follows from the various 
Definitions that there exists a natural identification $H_{i,n}= \widetilde \Pi_i[n]^{(\ab)}$ and a commutative diagram 
$$
\CD
1 @>>> \Delta_{i+1} @>>> \Pi_{i+1} @>>> G_k @>>> 1\\
@. @VVV @VVV   @|\\
1 @>>> \widetilde \Delta _{i,n}  @>>>\widetilde \Pi_{i,n} @>>> G_k @>>> 1\\
\endCD
$$
where the vertical maps are surjective. 
From this one deduces; since the section $s_{i,n}$ is induced by the section $s_i$ and $\Pi_{i,n}$ is a quotient of $\Pi_i$ [cf. diagram (2.8)],
the existence of a commutative diagram 

$$
\CD
1 @>>>  \Delta^{i+1}   @>>>  \Cal H_i@>>>  G_k  @>>> 1\\
@.    @VVV     @VVV     @| \\
1 @>>>  \widetilde \Delta _i[n]^{\ab}    @>>>   \widetilde \Pi_i[n]^{(\ab)} @>>>  G_k  @>>> 1\\
\endCD
$$
where the left and middle vertical map are surjective. Taking projective limits, and using the fact $\Delta ^{i+1}  \isom \underset{n\ge 1} \to{\varprojlim}\ \widetilde \Delta_i[n] ^{\ab}$,
it follows $\Cal H_i\isom \underset{n\ge 1} \to{\varprojlim}\  \widetilde \Pi_i[n]^{(\ab)}$. The case $(i,r)=(1,1)$ is treated in an entirely similar way as above.
\qed
\enddemo

Next, assume that the section $s_i:G_k\to \Pi_i$ lifts to a section $s_{i+1}:G_k\to \Pi_{i+1}$
of the natural projection $\Pi_{i+1}\twoheadrightarrow G_k$, i.e., there exists a section $s_{i+1}$ and a commutative diagram
$$
\CD
G_k  @>{s_{i+1}}>> \Pi_{i+1}\\
@|     @VVV\\
G_k  @>{s_i}>> \Pi _i\\
\endCD
$$
where the right vertical map is as in diagram (2.4). Recall the exact sequence (2.2)
$$1\to \Delta ^{i+1}\to \Delta _{i+1}\to \Delta _i\to 1.$$ 
We view $\Delta _{i+1}$
(hence also $\Delta _i$ and $\Delta ^{i+1}$) as a $G_k$-group via the action of $s_{i+1}(G_k)$ 
by conjugation. Thus, this $G_k$-group structure on $\Delta _i$ is the one induced by the section $s_i$.
The above sequence is an exact sequence of $G_k$-groups.

\proclaim {Lemma 2.3} Assume that $k$ is a {\bf number field} or a {\bf $p$-adic local field} 
(i.e., a finite extension of $\Bbb Q_p$). Then $H^0(G_k,\Delta ^{i+1})=\{0\}$, $\forall i\ge 0$.
\endproclaim

\demo{Proof} First, consider the case $i=0$ ($\Delta _{0}=\{1\}$, $\Delta^1=\Delta_1$). We have exact sequences of $G_k$-modules
$1\to I_X\to \Delta_1\to \Delta_1^{\et}\to 1$, and $1\to \hat \Bbb Z(1)\to \prod_{i=1}^m \Ind_{k(P_i)}^k\hat \Bbb Z(1)\to I_X\to 1$
(cf. discussion after Lemma 2.1). The group $H^0(G_k,\Delta_1 ^{\et})=\{0\}$ vanishes (see proof of Lemma 1.3 in [Sa\"\i di]). Further, we have an exact sequence
$0\to H^0(G_k, \hat \Bbb Z(1))\to H^0(G_k, \prod_{i=1}^m \Ind_{k(P_i)}^k \hat \Bbb Z(1))\to H^0(G_k,I_X)\to 0$ as follows from the fact that the map
$H^1(G_k, \hat \Bbb Z(1))\to H^1(G_k, \prod_{i=1}^m \Ind_{k(P_i)}^k\hat \Bbb Z(1))$ is injective. Indeed, this latter map is identified (via Kummer
theory) with the natural map $(k^{\times})^{\wedge}\to \prod_{i=1}^m (k(P_i)^{\times})^{\wedge}$; where for a field $\ell$
we write $(\ell^{\times})^{\wedge}\defeq \underset {j\ge 0} \to \varprojlim \ \ell^{\times}/(\ell^{\times})^j$, and we use Shapiro's Lemma
to identify $H^1(G_k,  \Ind_{k(P_i)}^k \hat \Bbb Z(1))$ with $H^1(G_{k(P_i)}, \hat \Bbb Z(1))$.
Now the map $(k^{\times})^{\wedge}\to \prod_{i=1}^m (k(P_i)^{\times})^{\wedge}$ is injective if $k$ is a number field or a $p$-adic local field.
This follows from the well-known structure of $k^{\times}$ if $k$ is a $p$-adic local field.
If $k$ is a number field this follows from the $p$-adic local field case and the fact that the map $(k^{\times})^{\wedge}\to \prod _{v\in \Primes_k} (k_v^{\times})^{\wedge}$
is injective (this follows from [Neukirch-Schmidt-Wingberg], (9.1.3)Theorem). 
Also $H^0(G_k, \hat \Bbb Z(1))=H^0(G_k, \prod_{i=1}^m \Ind_{k(P_i)}^k \hat \Bbb Z(1))=\{0\}$, which implies
$H^0(G_k,I_X)=\{0\}$, and hence $H^0(G_k,\Delta_1)=\{0\}$.

Next, assume $i\ge 1$. If $(i,r)\neq (1,1)$, the group $H^0(G_k,\Delta ^{i+1})$ is naturally identified with $\underset {n\ge 1}\to \varprojlim H^0(G_k, \widetilde \Delta_i[n]^{\ab})$
(cf. Proposition 2.2) where the $G_k$-module $\widetilde \Delta_i[n]^{\ab}$ is isomorphic to the Tate module of the jacobian of the curve $Y_{i,n}$.
Further, $H^0(G_k,\widetilde \Delta_i[n]^{\ab})=\{0\}$ (see proof of Lemma 1.3 in [Sa\"\i di]). Thus, $H^0(G_k,\Delta ^{i+1})=\{0\}$ in this case.
If $(i,r)=(1,1)$, the group $H^0(G_k,\Delta ^{i+1})$ in this case is naturally identified with $\underset {n\ge 1}\to \varprojlim H^0(G_k, \widetilde \Delta_i[n]^{\ab})$
(cf. Proposition 2.2) where the $G_k$-module $\widetilde \Delta_i[n]^{\ab}$ is by definition $\pi_1((X_{i,n})_{\bar k},\overline \eta_{i,n})^{\ab}$. The vanishing of 
$H^0(G_k,\Delta ^{i+1})$ follows in this case from the case $i=0$ discussed above.
\qed
\enddemo

\proclaim{Lemma 2.4} Assume that $k$ is a {\bf number field} or a {\bf $p$-adic local field} (i.e., a finite extension of $\Bbb Q_p$). Then $H^0(G_k,\Delta_i)=\{1\}$, $\forall i\ge 0$.
\endproclaim

\demo{Proof} By induction on $i$, using Lemma 2.3, and the fact that we have an exact sequence of groups
$1\to H^0(G_k,\Delta ^{i+1})\to H^0(G_k,\Delta _{i+1})\to H^0(G_k,\Delta _i)$.
\qed
\enddemo

\proclaim {Lemma 2.5} Assume that $k$ is a {\bf number field} or a {\bf $p$-adic local field} (i.e., a finite extension of $\Bbb Q_p$). Then the natural map
$H^1(G_k,\Delta ^{i+1})\to H^1(G_k,\Delta _{i+1})$ of pointed sets is {\bf injective}, $\forall i\ge 0$.
\endproclaim

\demo{Proof} There exists an exact sequence of pointed sets 
$$H^0(G_k,\Delta _i)\to H^1(G_k,\Delta ^{i+1})\to H^1(G_k,\Delta _{i+1})$$
(cf. [Serre], I, $\S5$, 5.5, Proposition 38). The proof follows from [Serre], I, $\S5$, 5.5, Proposition 39 (ii), and
the fact that $H^0(G_k,\Delta_i)$ is trivial (cf. Lemma 2.4).
\qed
\enddemo

\subhead
\S3. Proof of Theorem B
\endsubhead
This section is devoted to the proof of Theorem B (cf. $\S0$).
We use the same notations as in Theorem B, $\S0$, and $\S2$. We assume further that the set $S\subset \Primes_k$ 
contains no real place.

Recall the commutative diagram of exact sequences of profinite groups 
$$
\CD
1@>>> \Delta \defeq \pi_1(X_{\bar k},\overline \eta) @>>> \Pi\defeq \pi_1(X, \eta) @>>> G_k @>>> 1\\
@.  @VVV   @VVV   @|\\
1@>>> \Delta ^{\sol}@>>>   \Pi^{(\sol)}@>>> G_k@>>> 1\\
\endCD
$$
[cf. diagram (2.5)] and the natural map (cf. $\S0$)
$$\prod _{v\in S} \Res _v^{\sol}: H^1(G_k,\Delta^{\sol})\to \prod _{v\in S} H^1(G_{k_v},\Delta ^{\sol}).$$ 
We will show this map is {\bf injective}. Recall the definition of the $i$-th step solvable characteristic quotient $\Delta _i$ of $\Delta$ 
[cf. the discussion before the exact sequence (2.2)].

\proclaim {Proposition 3.1}
The natural map 
$$\prod _{v\in S} \Res _{v}^i: H^1(G_k,\Delta_i)\to \prod _{v\in S} H^1(G_{k_v},\Delta _i)$$ 
is {\bf injective} for $i\ge 1$.
\endproclaim

We will prove Proposition 3.1 by an induction argument on $i\ge 1$. The case $i=1$ is treated first in the following Proposition.

\proclaim {Proposition 3.2}
The natural map 
$$\prod _{v\in S} \Res _{v}^1: H^1(G_k,\Delta_1)\to \prod _{v\in S} H^1(G_{k_v},\Delta _1)$$ 
is {\bf injective}.
\endproclaim

\demo{Proof of Proposition 3.2}  Recall the exact sequence $0\to I_X\to \Delta_1\to \Delta_1^{\et}\to 0$ (cf. the discussion after Lemma 2.1).
We have a commutative diagram of exact sequences
$$
\CD
0 @>>> H^1(G_k,I_X) @>>> H^1(G_k,\Delta_1)@>>> H^1(G_k,\Delta_1^{\et}) \\
@.   @V\prod _{v\in S} \Res _{v}^1VV         @V\prod _{v\in S} \Res _{v}^1VV        @V\prod _{v\in S} \Res _{v}^{1,\et}VV\\
0 @>>> \prod _{v\in S} H^1(G_{k_v},I_X) @>>> \prod _{v\in S} H^1(G_{k_v},\Delta_1)@>>> \prod _{v\in S} H^1(G_{k_v},\Delta_1^{\et}) \\
\endCD
$$
where the horizontal sequences arise from the cohomology exact sequences associated to the exact sequence $0\to I_X\to \Delta_1\to \Delta_1^{\et}\to 0$
of $G_k$ (and $G_{k_v}$)-modules, the vertical maps are restriction maps, and the injectivity on the left in the horizontal sequences follows from Lemma 1.4 in [Sa\"\i di].
The injectivity of the left and right vertical maps in the above diagram would imply the injectivity of the middle vertical map.
The map $\prod _{v\in S} \Res _{v}^{1,\et}:H^1(G_k,\Delta_1^{\et}) \to \prod _{v\in S} H^1(G_{k_v},\Delta_1^{\et})$ is injective by Proposition 2.2 in [Sa\"\i di].
We shall prove the injectivity of the left vertical map $\prod _{v\in S} \Res _{v}^1:H^1(G_k,I_X) \to \prod _{v\in S} H^1(G_{k_v},I_X)$. 
To this end it suffices to prove, for a prime number $p$,  the injectivity of the map $\prod _{v\in S} \Res _{v}^1(p):H^1(G_k,I_X(p)) \to \prod _{v\in S} H^1(G_{k_v},I_X(p))$
where $I_X(p)$ is the $p$-primary part of $I_X$.
This follows from results of Serre
in [Serre1] and [Serre2]. 
For a profinite group $G$ and a continuous $G$-module $N$ write $H^1_{*}(G,N)\defeq \text {Ker} [H^1(G,N)@>\res>> \prod_C H^1(C,N)]$ where the product is over all pro-cyclic subgroups of $G$.
Recall the exact sequence $0\to \Bbb Z_p(1)@>>> \prod_{i=1}^n \Ind_{k(P_i)}^k \Bbb Z_p(1)\to I_X(p)\to 0$ of $G_k$-modules (cf. discussion after Lemma 2.1).
The $\Bbb Z_p$-module $I_X(p)$ is free of finite rank and is a $G_k$-module. Write $G\defeq \text {Image} [G_k\to \Aut(I_X(p))]$. Thus, $G$ is a $p$-adic Lie group. Write $\Cal G$
for the $p$-adic Lie algebra of $G$. The inflation map $H^1_{*}(G,I_X(p))\to H^1_{*}(G_k,I_X(p))$ is an isomorphism (cf. [Serre1], Proposition 6).
If $v\in \Primes _k$ is a prime of $k$ of residue characteristic $\neq p$, and $\text {Fr}_v\in G$ is a Frobenius at $v$, then its eigenvalues on $I_X(p)$ have complex absolute value $Nv$;
where $Nv$ is the cardinality of the residue field $\kappa(v)$ at $v$, and thus satisfy the condition in Lemme 2 of [Serre2].
Hence the Lie algebra $\Cal G$ of $G$ satisfies the hypothesis of Th\'eor\`eme 1 in [Serre2] and we have $H^1(\Cal G,I_X(p)\otimes_{\Bbb Z_p} \Bbb Q_p)=0$. This implies that
$H^1(G,I_X(p))$ is finite. Indeed, $H^1(G,I_X(p))$ is a finitely generated $\Bbb Z_p$-module and $H^1(G,I_X(p))\otimes_{\Bbb Z_p} \Bbb Q_p$ injects in
$H^1(\Cal G,I_X(p)\otimes_{\Bbb Z_p} \Bbb Q_p)=0$ (cf. [Serre2] Corollaire to Lemme 3). Moreover, since $(I_X(p)\otimes_{\Bbb Z_p}\Bbb Q_p)^{\text {Fr}_v}=0$ we have
$H^1_{*}(G,I_X(p))=0$ by [Serre1] Th\'eor\`eme 1. Now $\text {Kernel} [H^1(G_k,I_X(p))\to \prod _{v\in S} H^1(G_{k_v},I_X(p))]$ is contained in 
$H^1_{*}(G,I_X(p))$ (cf. Propositions 7 and 8 in [Serre1]), thus is trivial as required. Note that in [Serre1] and [Serre2] the set $S$ is the complement of a finite set 
of primes but the same proof works in the case where $S$ has density $1$ (see also [Jannsen], the proof of Theorem 3, and the references therein).
\qed
\enddemo

\demo {Proof of Proposition 3.1} The proof is similar to the proof of Proposition 2.1 in [Sa\"\i di] using Proposition 3.2 and Proposition 2.2.
Fix an integer $i\ge 1$.  Consider the following commutative diagram of maps of pointed cohomology sets
$$
\CD
1 @>>> H^1(G_k,\Delta ^{i+1}) @>>> H^1(G_k,\Delta _{i+1}) @>>> H^1(G_k,\Delta _{i})  \\
@.   @VVV   @V{\prod _{v\in S} \Res _{v}^{i+1}}VV   @V{\prod _{v\in S} \Res _{v}^i}  VV \\
1 @>>> \prod _{v\in S} H^1(G_{k_v},\Delta ^{i+1}) @>>> \prod _{v\in S} H^1(G_{k_v},\Delta _{i+1})@>>> \prod _{v\in S} H^1(G_{k_v},\Delta _{i})\\
\endCD
$$
where the horizontal sequences are exact (cf. Lemma 2.5) and the vertical maps are the natural restriction maps.
We assume by induction hypothesis that the right vertical map
$\prod _{v\in S} \Res _{v}^i: H^1(G_k,\Delta_i)\to \prod _{v\in S} H^1(G_{k_v},\Delta _i)$
is {\bf injective}. We will show that the middle vertical map
$\prod _{v\in S} \Res _{v}^{i+1}: H^1(G_k,\Delta_{i+1})\to \prod _{v\in S} H^1(G_{k_v},\Delta _{i+1})$ 
is {\bf injective}.
Let 
$$[\rho], [\tau]\in H^1(G_k,\Delta_{i+1})$$
be two cohomology classes such that
$$\prod _{v\in S} \Res _{v}^{i+1}([\tau])=\prod _{v\in S} \Res _{v}^{i+1}([\rho]).$$
Write $s_j:G_k\to \Pi_{j}$ for the section of the projection $\Pi_{j}\twoheadrightarrow G_k$ induced by the section $s:G_k\to \Pi$, for $j\ge 1$.
We will show $[\tau]=[\rho]$. We can (without loss of generality), and will, assume that $[\tau]=[s_{i+1}]=1$ is the distinguished element of $H^1(G_k,\Delta _{i+1})$.
The classes $[\rho]$ and $[\tau]$ map to the classes $[\rho_1]$ and $[\tau_1]=[s_i]=1$ in $H^1(G_k,\Delta _i)$, respectively.
In particular, we have the equality
$$\prod _{v\in S} \Res _{v}^i([\tau_1])=\prod _{v\in S} \Res _{v}^i([\rho_1]),$$
hence $[\tau_1]=[\rho_1]=1$ since $\prod _{v\in S} \Res _{v}^i$ is injective by the induction hypothesis. 
Thus, there exist classes $[\tilde \tau]=1$ and $[\tilde \rho]$ in $H^1(G_k,\Delta ^{i+1})$ which map to $[\tau]$ and $[\rho]$, respectively
in $H^1(G_k,\Delta _{i+1})$ (cf. exactness of horizontal sequences in the above diagram).

Next, and in order to show that $[\rho]=[\tau]$ in $H^1(G_k,\Delta_{i+1})$
it suffices to show $[\tilde \rho]=[\tilde \tau]$ in $H^1(G_k,\Delta^{i+1})$.
Note that the assumption $\Res_v^{i+1}([\tau])= \Res_v^{i+1}([\rho])$
implies that  $\Res_v^{i+1} ([\tilde \tau])= \Res_v^{i+1}([\tilde \rho])$ in $H^1(G_{k_v},\Delta ^{i+1})$, for each place $v\in S$,
as follows from the injectivity of the maps $H^1(G_{k_v},\Delta^{i+1})\to H^1(G_{k_v},\Delta_{i+1})$
(cf. Lemma 2.5). 
We have a commutative diagram of group homomorphisms
$$
\CD
H^1(G_k,\Delta^{i+1})  @>>> \underset{n\ge 1} \to{\varprojlim} H^1(G_k, \widetilde \Delta _i[n]^{\ab})\\
@V{\prod _{v\in S} \Res _{v}^{i+1}}VV                      @VVV\\
\prod _{v\in S}H^1(G_{k_v},\Delta ^{i+1})        @>>>     \prod _{v\in S}  \lgroup \underset{n\ge 1} \to{\varprojlim}  H^1(G_{k_v},\widetilde \Delta _i[n]^{\ab})\rgroup\\
\endCD
$$
where the horizontal maps are induced by the identification 
$\Cal H_i\isom \underset{n\ge 1} \to{\varprojlim}\  \widetilde \Pi_i[n]^{(\ab)}$ (cf. Proposition 2.2) and the vertical maps are restriction maps.
The identification $\Cal H_i\isom \underset{n\ge 1} \to{\varprojlim}\  \widetilde \Pi_i[n]^{(\ab)}$
induces isomorphisms
$H^1(G_k,\Delta^{i+1})\isom \underset{n\ge 1} \to{\varprojlim} H^1(G_k, \widetilde \Delta _i[n]^{\ab})$,
and $H^1(G_{k_v},\Delta^{i+1})\isom \underset{n\ge 1} \to{\varprojlim} H^1(G_{k_v}, \widetilde \Delta _i[n]^{\ab})$
for each $v\in S$, where the transition homomorphisms in the projective limit are induced by the natural
$G_k$ (resp. $G_{k_v}$)-homomorphisms $\widetilde \Delta _{i}[n+1]^{\ab}\to \widetilde \Delta _{i}[n]^{\ab}$ (cf. [Neukirch-Schmidt-Wingberg], Chapter II, (2.3.5)Corollary). 
The right vertical map in the above diagram is injective in the case $(i,r)\neq (1,1)$ by Proposition 2.2 in [Sa\"\i di]
[recall in this case that the $G_k$-module $\widetilde \Delta_i[n] ^{\ab}$ is isomorphic to the Tate module of the jacobian of the curve $Y_{i,n}$ defined in $\S2$].
In case $(i,r)=(1,1)$ the injectivity of the right vertical map in the above diagram follows from Proposition 3.2 [recall; with the notations in $\S2$, that in this case the $G_k$-module 
$\widetilde \Delta_i[n] ^{\ab}$ is $\pi_1((X_{i,n})_{\bar k},\overline \eta_{i,n})^{\ab}$].

Thus, the left vertical map
$H^1(G_k,\Delta ^{i+1})\to \prod _{v\in S}H^1(G_{k_v},\Delta ^{i+1})$
is {\bf injective}. In particular, $[\tilde \rho]=[\tilde \tau]$ and $[\rho]=[\tau]$.
This finishes the proof of Proposition 3.1.
\qed
\enddemo

Finally, going back to the proof of Theorem B, let 
$$[\alpha], [\beta]\in H^1(G_k,\Delta^{\sol})$$
be two cohomology classes such that
$$\prod _{v\in S} \Res _{v}^{\sol}([\alpha])=\prod _{v\in S} \Res _{v}^{\sol}([\beta]).$$
One shows $[\alpha]=[\beta]$ by similar arguments to the ones used to conclude the proof of Theorem B in [Sa\"\i di] (cf. loc. cit. discussion after the proof of Lemma 2.1),
using Proposition 2.1 and Lemma 2.3 in loc. cit..

This finishes the Proof of Theorem B.
\qed

$$\text{References.}$$

\noindent
[Grothendieck] Grothendieck, A., Rev\^etements \'etales et groupe fondamental, Lecture 
Notes in Math. 224, Springer, Heidelberg, 1971.

\noindent
[Jannsen] Jannsen, U., On the $l$-adic cohomology of varieties over number fields and Galois cohomology. Galois groups over $\Bbb Q$, edited by
Y. Ihara, K. Riber, and J.-P. Serre, Proceedings of a Workshop Held
March 23-27, 1987, Mathematical Sciences Research Institute Publications 16, Springer-Verlag.



\noindent
[Neukirch-Schmidt-Wingberg] Neukirch, J., Schmidt, A., Wingberg, K., Cohomology of number fields, first edition,
Springer, Grundlehren der mathematischen Wissenschaften Bd. 323, 2000.


\noindent
[Sa\"\i di] Sa\"\i di, M., A local-global principle for torsors under geometric prosolvable fundamental groups, 
Manuscripta Mathematica 145, no. 1-2 (2014), 163-174.

\noindent
[Sa\"\i di1] Sa\"\i di, M., Arithmetic of $p$-adic curves and sections of geometrically abelian fundamental groups. Math. Z. 297 (2021), no. 3-4, 1191-1203.




\noindent
[Serre] Serre, J.-P., Galois cohomology, Springer-Verlag Berlin Heidelberg, 1997.

\noindent
[Serre1] Serre, J.-P., Sur Ies groupes de congruence des vari\'et\'es ab\'eliennes I, Izv. Akad.
Nauk SSSR, S\'er. Mat. 28 (1964), 3-20.

\noindent
[Serre2] Serre, J.-P., Sur Ies groupes de congruences des vari\'et\'es ab\'eliennes II, Izv. Akad.
Nauk SSSR, S\'er. Mat. 35 (1971), 731-737.

\bigskip

\noindent
Mohamed Sa\"\i di

\noindent
College of Engineering, Mathematics, and Physical Sciences

\noindent
University of Exeter

\noindent
Harrison Building

\noindent
North Park Road

\noindent
EXETER EX4 4QF 

\noindent
United Kingdom

\noindent
M.Saidi\@exeter.ac.uk

\end
\enddocument